\documentclass[12pt]{amsart}
\usepackage{amssymb}
\usepackage[english]{babel}
\usepackage{enumerate}
\usepackage{dsfont}

\def\C{{\mathbb C}}

\def\D{{\mathbb D}}
\def\L{{\mathcal L}}

\newcommand{\norm}[1]{\left\Vert#1\right\Vert}
\DeclareMathOperator{\dummydist}{dist}
\newcommand{\dist}[1]{\dummydist ( #1 )}

\DeclareMathOperator{\cabv}{cabv}

\newtheorem{theorem}{Theorem}[section]

\newtheorem{lemma}[theorem]{Lemma}
\newtheorem{corollary}[theorem]{Corollary}

\newtheorem{assumption}{Assumption}

\theoremstyle{definition}

\newtheorem{remark}[theorem]{Remark}

\theoremstyle{empty}

\begin{document}

\title{On M-ideals and \MakeLowercase{o}--O type spaces}

\author{Karl-Mikael Perfekt}

\address{%
  Department of Mathematical Sciences \\
  Norwegian University of Science and Technology \\
  7491 Trondheim
}

\thanks{2010 Mathematics Subject Classification: 46B04, 46B10, 46B25. }

\begin{abstract}
We consider pairs of Banach spaces $(M_0, M)$ such that $M_0$ is defined in terms of a little-$o$ condition, and $M$ is defined by the corresponding big-$O$ condition. The construction is general and pairs include
function spaces of vanishing and bounded mean oscillation, vanishing weighted and weighted spaces of functions or their derivatives, M\"obius invariant spaces of analytic functions, Lipschitz-H\"older spaces, etc.
It has previously been shown that the bidual $M_0^{**}$ of $M_0$ is isometrically isomorphic with $M$. The main result of this paper is that $M_0$ is an M-ideal in $M$. This has several useful consequences:
$M_0$ has Pe\l cz\'{y}nskis properties (u) and (V), $M_0$ is proximinal in $M$, and $M_0^*$ is a strongly unique predual of $M$, while $M_0$ itself never is a strongly unique predual.
\end{abstract}

\maketitle 

\section{Introduction}
The aim of this work is to show that Banach spaces whose definitions are given in terms of little-$o$ conditions are M-embedded. That is, to show that they are M-ideals in their bidual spaces, the latter spaces which may be canonically identified with the Banach spaces defined by the corresponding big-$O$ conditions. We will treat a large class of spaces, our main result yielding that a vast array of classical spaces studied in analysis in fact turn out to be examples of M-ideals: spaces of vanishing mean oscillation, vanishing weighted spaces of continuous, harmonic, or analytic functions or their derivatives, the little versions of general M\"obius invariant spaces of analytic functions, Lipschitz-H\"older spaces, and many more. 

The notion of the M-ideal, as a Banach space analogue of a two-sided ideal in a $C^*$-algebra, was born in Alfsen's and Effros' influential paper \cite{Alfsen72}. As for M-embedded spaces, their systematic study was initiated by Harmand and Lima \cite{Harm84}. We refer to the comprehensive monograph of Harmand, D. Werner, and W. Werner \cite{Harm93}, not only for further notes on the literature, but also for an excellent presentation of the available theory of M-ideals. 

From the point of view of this paper, showing that a Banach space is M-embedded carries the benefit of the immediate application of the rich theory associated with M-ideals. For instance, M-embedded spaces have Pe\l cz\'{y}nskis properties (u) and (V), which the author in \cite{Perf15} utilized to characterize all weakly compact operators acting on spaces defined by
little-$o$ conditions. Further examples of the strong geometric results available for an M-embedded Banach space $Z$ are given by the facts that $Z$ is always proximinal in $Z^{**}$ and that $Z^*$ is the strongly unique predual of $Z^{**}$. We shall return to these applications later in this section, as corollaries of the main result. 

The present work is motivated by the fact that known examples of non-reflexive M-embedded Banach spaces $Z$ often have the character of a little space -- "vanishing at infinity" in some sense, if one permits the use of vague terminology -- while the space $Z^{**}$ acts as the corresponding big space. This is of course exhibited by the archetypal M-embedded space, namely, the sequence space $c_0$; $c_0$ is an M-ideal in $c_0^{**} = \ell^\infty$. To observe similar behavior of many other concrete examples of M-embedded spaces, we refer for example to (\cite{Harm93}, III.1) \cite{Kaminska04}, \cite{Lueck80}, or \cite{Werner92}.
 
The goal of this article is therefore, in a sense, to formalize the intuition presented in the previous paragraph.   In \cite{Perf13}, the author considered a general construction of pairs of Banach spaces $(M_0, M)$ -- a little space $M_0$ defined by a little-$o$ condition, and a big space $M$ defined by the corresponding big-$O$ condition. One of the main results of the aforementioned paper is that $M_0^{**} \simeq M$ in a canonical way. The main theorem of the present work states that $M_0$ is in fact an M-ideal in $M$. This gives a new range of concrete examples of M-embedded spaces taken from harmonic and complex analysis, since examples of pairs $(M_0, M)$ include vanishing and bounded mean oscillation in one and more variables, general M\"{o}bius invariant spaces of holomorphic function, and Lipschitz-H\"older spaces. Note that these spaces are all considered with their instrinsic norms. We will in a moment define the spaces $M_0$ and $M$, but we refer to \cite{Perf13} for a detailed treatment of the realization of these examples within the framework. 

The definition of $(M_0, M)$ relies on several auxiliary objects, which we now fix. Let $X$ and $Y$ be two Banach spaces, where $X$ is separable and reflexive. The norm of $M$ will be determined through a collection $\L \subset B(X, Y)$ of bounded linear operators $L : X \to Y$. By equipping $\L$ with a topology $\tau$ we are able to give meaning to the statement that elements of $M_0$ vanish at infinity. The topological space $(\L, \tau)$ should be Hausdorff, $\sigma$-compact, and locally compact, and for every $x \in X$ the map $L \mapsto Lx$ should act continuously from $(\L, \tau)$ to $Y$. The limit $L \to \infty$ is now given the standard meaning of $L$ escaping all compact sets of $(\L, \tau)$, or equivalently that $L$ tends to $\infty$ in the one-point compactification $\alpha \L = \L \cup \{ \infty \}$ of $\L$.

The spaces $M$ and $M_0$ are defined by
\begin{equation} \label{eq:mdef}
M(X,\L) = \left\{ x \in X \, : \, \sup_{L \in \L} \norm{Lx}_Y < \infty \right\}
\end{equation}
and
\begin{equation} \label{eq:m0def}
M_0(X,\L) = \left\{ x \in M(X, \L) \, : \, \varlimsup_{\L \ni L\to\infty} \norm{Lx}_Y = 0 \right\}.
\end{equation}
We assume that $M(X, \L)$ is dense in $X$ under the $X$-norm, and that $M(X, \L)$ is a Banach space continuously contained in $X$ under the norm
\begin{equation*}
\norm{x}_M = \sup_{L \in \L} \norm{Lx}_Y.
\end{equation*}
To ask the question whether $M_0$ is M-embedded by being an M-ideal in $M$, we must first isometrically identify the bidual space $M_0^{**}$ with $M$. In \cite{Perf13} it was shown that $M_0^{**}$ is canonically isometrically isomorphic with $M$ (see Theorem \ref{thm:perf13}) if and only if we have the following approximation property, which we refer to as Assumption \ref{as2}. In the sequel we always assume that Assumption \ref{as2} holds.

\begin{assumption} \label{as2} For every $x\in M(X,\L)$ there is a bounded sequence $(x_n)_{n=1}^\infty$ in $M_0(X,\L)$ such that $x_n $ converges weakly to $x$ in $X$ and $\sup_n \norm{x_n}_{M(X,\L)} \leq \norm{x}_{M(X,\L)}$. 
\end{assumption}

We are now in a position to state the main theorem.
\begin{theorem} \label{thm:main}
Suppose that Assumption \ref{as2} holds. Then $M_0(X, \L)$ is an $M$-embedded Banach space. That is, it is an M-ideal in $M_0(X, \L)^{**} \simeq M(X, \L)$. 
\end{theorem}
As mentioned previously, Theorem \ref{thm:main} has a number of immediate corollaries. In \cite{Perf13} the distance between an element $x \in M$ and the space $M_0$ was computed. Since $M$-embedded spaces are always proximinal in their biduals \cite{Alfsen72}, \cite{Ando73} (the distance between an element of the bidual and the space has a least minimizer), we obtain in conjunction with the distance calculation the following result.
\begin{corollary}
For every $x \in M(X, \L)$ it holds that 
\begin{equation*} 
\dist{x, M_0(X, \L)}_{M(X,\L)} = \min_{x_0 \in M_0} \|x-x_0\|_M = \varlimsup_{\L \ni L\to\infty} \norm{Lx}_Y.
\end{equation*}
\end{corollary}
A Banach space $Z$ is said to be the strongly unique predual of $Z^*$ if every isometric isomorphism from $Z^*$ onto $W^*$, $W$ a Banach space, is the adjoint of an isometric isomorphism of $W$ onto $Z$. From Proposition 2.10 of (\cite{Harm93}, p. 122) we obtain the following corollary. The reflexive case $M_0 = M$ has to be excluded.
\begin{corollary} \label{cor:uniquepredual}
Suppose that $M_0(X, \L) \neq M(X, \L)$. Then
\begin{enumerate}
\item $M_0(X, \L)^*$ is the strongly unique predual of $M(X, \L)$.
\item $M_0(X, \L)$ is never a strongly unique predual.
\end{enumerate}
\end{corollary}
\begin{remark}
Part (1) of Corollary \ref{cor:uniquepredual} was previously shown, with a different proof, in \cite{Perf13}.
\end{remark}
Theorem \ref{thm:perf13} implies that $M_0^*$ is separable, hence also that $M_0$ always is a separable space. Godefroy and Li \cite{Godefroy90} proved that a separable M-embedded space is an $\L^\infty$ space (see for instance \cite{Linden69}) if and only if it is isomorphic to $c_0$. 
\begin{corollary}
If $M_0(X, \L)$ is an $\L^\infty$ space, then $M_0(X, \L)$ is isomorphic to $c_0$ and $M(X,\L)$ is isomorphic to $\ell^\infty$.
\end{corollary}
\begin{remark}
Let $\D$ denote the unit disk in the complex plane $\C$, and denote by $v : [0,1] \to [0,\infty]$ a continuous, decreasing weight function such that $v(1) = 0$. The vanishing weighted space of holomorphic functions
$$(Hv)_0 = \{f : \D \to \C \textrm{ holomorphic} \, : \, \varlimsup_{|z| \to 1} |f(z)|v(|z|) = 0 \}$$
is a basic example of a space of the form $M_0$. Lusky \cite{Lusky06} has completely characterized the weights $v$ for which $(Hv)_0$ is isomorphic to $c_0$.
\end{remark}
For the final corollary, we note that M-embedded spaces possess Pe\l cz\'{y}nskis properties (V) \cite{Godefroy89} and (u) \cite{Godefroy892}. We hence obtain the following, which is restatement of the fact that $M_0$ has property (V) (see (\cite{Harm93}, p. 128)).
\begin{corollary} \label{cor:weakcpct}
If $Z$ is a Banach space and $T : M_0(X,\L) \to Z$ is a bounded operator, then $T$ is weakly compact if and only if there does not exist a subspace $F \subset M_0(X, \L)$ isomorphic to $c_0$ such that $T|_F$ is an isomorphism. 
\end{corollary}
\begin{remark}
Several recent papers \cite{Cont14}, \cite{Lait11}, \cite{Lait13}, \cite{Lefevre10} have made use of the construction of $c_0$-subspaces to characterize the compactness of composition and integration operators acting on spaces of analytic functions of $M_0$ type. These concrete operators all exhibit the behavior of being compact precisely when weakly compact. This is investigated further in \cite{Perf15}.
\end{remark}
The remainder of this paper is organized as follows. Section 2 discusses preliminaries of the spaces $M_0$ and $M$, M-ideals, and some vector-valued integration theory. In Section 3 the main result is proven.
\section{Definitions and preliminaries} 
\subsection{The spaces $M_0$ and $M$} $M(X, \L)$ and $M_0(X, \L)$ were previously defined in \eqref{eq:mdef} and \eqref{eq:m0def}, but we now recall the precise formulation of the fact that $M_0^{**} \simeq M$ isometrically. For the statement, note that $M_0$ can be considered a closed subspace of both $M_0^{**}$ and $M$.
\begin{theorem}[\cite{Perf13}] \label{thm:perf13}
Suppose that Assumption \ref{as2} holds. Then $X^*$ is continuously contained and dense in $M_0(X, \L)^*$. Denoting by $$I \colon  M_0(X, \L)^{**} \to X$$ the adjoint of the inclusion map $J \colon  X^* \to M_0(X,\L)^*$, the operator $I$ is an isometric isomorphism of $M_0(X, \L)^{**}$ onto $M(X,\L)$ which acts as the identity on $M_0(X,\L)$.  
\end{theorem}

\subsection{M-ideals}
Suppose that $Z$ is a Banach space. A (closed) subspace $J \subset Z$ is called an M-ideal if the annihilator $J^\perp \subset Z^*$ is the range of an L-projection -- a projection $L : Z^* \to Z^*$ such that
\begin{equation*}
\|z^*\| = \|Lz^*\| + \|z^* - Lz^*\|, \quad \forall z^* \in Z^*.
\end{equation*}
An $M$-embedded space $Z$ is a Banach space which is an $M$-ideal when considered as a subspace of its bidual $Z^{**}$. Note that there is always a canonical projection $\pi : Z^{***} \to Z^{*}$ with range $Z^*$ and kernel $Z^\perp \subset Z^{***}$,
\begin{equation*}
(\pi z^{***})(z) = z^{***}(z), \quad z^{***} \in Z^{***}, \, z \in Z. 
\end{equation*}
Here and in the sequel we freely consider any Banach space to be a subspace of its bidual without special notation. It is a basic fact (\cite{Harm93}, p. 102) that $Z$ is an $M$-ideal in $Z^{**}$ if and only if the canonical projection $\pi$ is an L-projection. Hence the fact that $Z$ is $M$-embedded is equivalently expressed by saying that the canonical decomposition $Z^{***} = Z^* \oplus Z^{\perp}$ induced by $\pi$ is an $\ell^1$-decomposition,
\begin{equation*} \label{decomp}
Z^{***} = Z^* \oplus_1 Z^{\perp}.
\end{equation*}

\subsection{Measure theory}
The proof of Theorem \ref{thm:main} relies on studying duality via the embedding $V : M(X, \L) \to C_b(\L, Y)$,
\begin{equation*}
(Vx)(L) = Lx, \quad x \in M, \, L \in \L.
\end{equation*}
Here  $C_b(\L, Y)$ denotes the space of bounded continuous $Y$-valued functions on $(\L, \tau)$, equipped with the supremum norm 
\begin{equation*}
\| T \|_{C_b} = \sup_{L \in \L} \| T(L) \|_Y, \quad T \in C_b(\L, Y).
\end{equation*}
Note that $V$ isometrically embeds $M(X, \L)$ into $C_b(\L, Y)$ and that it similarly embeds $M_0(X, \L)$ into the space $C_0(\L, Y)$ of continuous functions vanishing at $\infty$.

We will require a few elements of $Y$-valued measure theory. We refer to \cite{dobrakov3}, \cite{Perf13}, and \cite{Zinger}. The space of countably additive $Y^*$-valued Baire measures of bounded variation is denoted by $\textrm{cabv}(\L, Y^*)$. It is equipped with the usual variation norm
\begin{equation*}
\norm{\mu}_{\textrm{cabv}} = \sup \sum \norm{\mu(\mathcal{E}_i)}_{Y^*} < \infty,
\end{equation*} 
where the supremum is taken over all pairwise disjoint partitions of $\L$ into sets $\mathcal{E}_i$. 

The reason for introducing $\textrm{cabv}(\L, Y^*)$ is of course the Riesz-Zinger theorem; $\textrm{cabv}(\L, Y^*)$ is isometrically isomorphic with the dual space $C_0(\L, Y)^*$ and we will freely identify the two. To be more precise about the identification, we introduce the pairing $\langle T, \mu \rangle$ between a function $T \in C_b(\L, Y)$ and a measure $\mu \in \textrm{cabv}(\L, Y^*)$,
\begin{equation} \label{eq:pairing}
\langle T, \mu \rangle = \int_\L T(L) \, d \mu(L).
\end{equation}

\begin{theorem}[\cite{dobrakov3}, \cite{Perf13}] \label{thm:zinger}
For every $\ell \in C_0(\L, Y)^*$ there is a unique measure $\mu \in \cabv(\L, Y^*)$ such that $\ell(T) = \langle T, \mu \rangle$ for all $T \in C_0(\L, Y)$. Conversely, each measure $\mu$ defines an element $\ell \in C_0(\L, Y)^*$ through \eqref{eq:pairing}, and $\norm{\ell}_{C_0^*} = \norm{\mu}_{\textrm{cabv}}$. 

Furthermore, each $T \in C_b(\L, Y)$ defines an element $k \in \cabv(\L, Y^*)^*$ by the formula $k(\mu) = \langle T, \mu \rangle$, and  $\|k\|_{\textrm{cabv}^*} = \|T\|_{C_b}$. The isometric embedding of $C_b(\L, Y)$ into $\cabv(\L, Y^*)^*$ given by $T \mapsto k$ extends the canonical embedding of $C_0(\L, Y)$ into $C_0(\L, Y)^{**}$.
\end{theorem}

\section{Proof of the main theorem} 
We begin by explaining the notation to be used in the proof of Theorem \ref{thm:main}. For $m \in M(X,\L)^*$,  $m \circ V^{-1}$ acts on $VM(X,\L) \subset C_b(\L, Y)$, which we as in Theorem \ref{thm:zinger} view as a subspace of $\cabv(\L,Y^*)^*$. By the Hahn-Banach theorem, $m \circ V^{-1}$ hence extends to a functional $\bar{m} \in \cabv(\L,Y^*)^{**}$ satisfying $\norm{\bar{m}} = \norm{m}$. Applying the canonical decomposition $Z^{***} = Z^* \oplus Z^{\perp}$ with $Z = C_0(\L, Y)$ we obtain
\begin{equation*}
\cabv(\L,Y^*)^{**} = \cabv(\L,Y^*) \oplus C_0(\L,Y)^\perp,
\end{equation*}
and we decompose $\bar{m}$ accordingly, 
$$\bar{m} = \bar{m}_{\omega^*} + \bar{m}_s, \quad \bar{m}_{\omega^*} \in \cabv(\L,Y^*), \; \bar{m}_s \in C_0(\L,Y)^\perp.$$

On the other hand, letting $I : M_0(X,\L)^{**} \to M(X,\L)$ be the isometric isomorphism of \ref{thm:perf13}, we obtain a second decomposition $m \circ I = (m \circ I)_{\omega^*} + (m \circ I)_s$ from
\begin{equation*} \label{mxdualdecomp} 
M(X,\L)^* \simeq M_0(X,\L)^{***} = M_0(X,\L)^* \oplus M_0(X,\L)^\perp.
\end{equation*}
Here $(m \circ I)_{\omega^*} \in M_0(X,\L)^*$ and $(m \circ I)_s \in M_0(X,\L)^\perp$.

Claim 3.5 of \cite{Perf13} amounts to the fact that the first decomposition is an extension of the second. We restate this here, as a lemma.
\begin{lemma} \label{lem:decomp}
In the above notation, we have $$\bar{m}_{\omega^*} \circ V \circ I = (m \circ I)_{\omega^*}$$ and 
$$\bar{m}_s \circ V \circ I = (m \circ I)_s,$$ 
as functionals on $M_0(X, \L)^{**}$. 
\end{lemma}

We are now prepared to prove the main theorem.

{
\renewcommand{\thetheorem}{\ref{thm:main}}
\begin{theorem}
  $M_0(X, \L)$ is an M-ideal in $M(X, \L)$. That is,
  \begin{equation*} 
  M(X,\L)^* \simeq M_0(X,\L)^{***} = M_0(X,\L)^* \oplus_1 M_0(X,\L)^\perp .
  \end{equation*}
\end{theorem}
\addtocounter{theorem}{-1}
}

\begin{proof}
Let $h \in M_0(X, \L)^{***}$ and define $m \in M(X,\L)^*$ by  $m = h \circ I^{-1}$. We employ the notation of this section, so that constructs involving $m$ are defined as above. Let $\mu \in \cabv(\L,Y^*)$ be the measure corresponding to $\bar{m}_{\omega^*}$, which in particular means that 
\begin{equation*}
\bar{m}_{\omega^*}(T) = \int_\L T(L) \, d \mu(L), \quad T \in C_b(\L,Y).
\end{equation*}
Denote by $\ell$ the restriction of $\bar{m}_s$ to $C_b(\L, Y)$, and let
\begin{equation*}
\tilde{m} = \bar{m}|_{C_b(\L, Y)} = \mu + \ell.
\end{equation*}
Here and in the remainder of the proof we understand $\mu$ as a functional on $C_b(\L, Y)$, as well as a measure in $\cabv(\L,Y^*)$, by slight abuse of notation which is justified in Theorem \ref{thm:zinger}.

Let $\mathcal{K}_1 \subset \mathcal{K}_2 \subset \cdots$ be an increasing sequence of compact Baire measurable subsets of $(\L, \tau)$ such that $\L = \bigcup_{n=1}^\infty \mathcal{K}_n$.  Denote, as before, by $\alpha \L = \L \cup \{\infty\}$ the one point compactification of $\L$. For each $n$, let $s_n \colon  \alpha L \to [0,1]$ be a continuous function such that $s_n^{-1}(1) \supset \mathcal{K}_n$ and $s_n(\infty) = 0$.

Now let $\mu_n = \mu|_{\mathcal{K}_n}$ be the restriction of the measure $\mu$ to $\mathcal{K}_n$, and consider the functional $\tilde{m}_n = \mu_n + \ell$ acting on $C_b(\L, Y)$. For fixed $n$, given $\varepsilon > 0$, let $S, T \in C_b(\L, Y)$ be such that
\begin{equation*}
\norm{S}_{C_b} = \norm{T}_{C_b} = 1,\quad \mu_n(S) > \norm{\mu_n}_{C_b^*} - \varepsilon, \quad \ell(T) > \norm{\ell}_{C_b^*} - \varepsilon.
\end{equation*}
Note that by construction we have
\begin{align*}
\tilde{m}_n(s_n S + (1-s_n)T) &= \mu_n(s_n S) + \ell((1-s_n) T) \\ &= \mu_n(S) + \ell(T) > \norm{\mu_n}_{C_b^*} + \norm{\ell}_{C_b^*} - 2\varepsilon .
\end{align*}
Observing that $\norm{ s_n S + (1-s_n)T }_{C_b} \leq 1$, we deduce in letting $\varepsilon \to 0$ that
\begin{equation*}
\norm{\tilde{m}_n}_{C_b^*} = \norm{\mu_n}_{C_b^*} + \norm{\ell}_{C_b^*}.
\end{equation*}
Since $\|\mu - \mu_n \|_{C_b^*} \leq \|\mu - \mu_n \|_{\textrm{cabv}} \to 0$,  we obtain by also letting $n \to \infty$ that
\begin{equation*}
\norm{\tilde{m}}_{C_b^*} = \norm{\mu}_{C_b^*} + \norm{\ell}_{C_b^*}.
\end{equation*}
Since $\norm{m}_{M^*} = \norm{\tilde{m}}_{C_b^*} = \norm{\bar{m}}_{\textrm{cabv}^{**}}$, we conclude that
\begin{equation*}
\norm{m}_{M^*} = \norm{\mu}_{C_b^*} + \norm{\ell}_{C_b^*}.
\end{equation*}
In view of Lemma \ref{lem:decomp} we get the following inequality:
\begin{align*}
\norm{h}_{M_0^{***}} &= \norm{m}_{M^*} \\ &= \norm{\mu}_{C_b^*} + \norm{\ell}_{C_b^*} \\ &= \norm{\bar{m}_{\omega*}}_{C_b^*} + \norm{\bar{m}_{s}}_{C_b^*} \\
&\geq \norm{\bar{m}_{\omega*}}_{(VM)^*} + \norm{\bar{m}_{s}}_{(VM)^*} \\ &= \norm{\bar{m}_{\omega*} \circ V}_{M^*} + \norm{\bar{m}_{s} \circ V}_{M^*} \\ 
&= \norm{\bar{m}_{\omega*} \circ V \circ I}_{M_0^{***}} + \norm{\bar{m}_{s} \circ V \circ I}_{M_0^{***}} \\ &= \norm{(m \circ I)_{\omega*}}_{M_0^{***}} + \norm{(m \circ I)_{s}}_{M_0^{***}} 
\\ &= \norm{h_{\omega*}}_{M_0^{***}} + \norm{h_{s}}_{M_0^{***}},
\end{align*}
where 
$$h_{\omega*} = (m \circ I)_{\omega*} \in M_0(X, \L)^*, \quad h_{s} = (m \circ I)_{s} \in M_0(X, \L)^\perp,$$ and $h = h_{\omega*} + h_s$. To finish the proof we only need to note the opposite inequality $$\norm{h}_{M_0^{***}} \leq \norm{h_{\omega*}}_{M_0^{***}} + \norm{h_{s}}_{M_0^{***}},$$ which is obvious.
\end{proof}

\end{document}